\documentclass[reqno,12pt]{amsart}

\usepackage{epsf}
\usepackage{graphics}
\usepackage{amssymb}
\usepackage{amsmath}

\date{}

\theoremstyle{plain}
\newtheorem{theorem}{Theorem}

\newtheorem{lemma}{Lemma}
\newtheorem{proposition}{Proposition}

\theoremstyle{definition}
\newtheorem*{definition}{Definition}

\theoremstyle{remark}

\newtheorem*{remark}{Remark}

\def\C{{\mathbb C}}

\def\N{{\mathbb N}}

\def\Z{{\mathbb Z}}

\def\Q{{\mathbb Q}}

\def\R{{\mathbb R}}

\title{Asymptotic link invariants for ergodic vector fields} 

\author{Sebastian Baader}

\begin{document}

\begin{abstract} We study the asymptotics of a family of link invariants on the orbits of a smooth volume-preserving ergodic vector field on a compact domain of the 3-space. These invariants, called linear saddle invariants, include many concordance invariants and generate an infinite-dimensional vector space of link invariants. In contrast, the vector space of asymptotic linear saddle invariants is 1-dimensional, generated by the asymptotic signature. We also relate the asymptotic slice genus to the asymptotic signature.
\end{abstract}

\maketitle

\section{Introduction}

A smooth vector field $X$ on a domain $G \subset \R^3$ produces a foliation with singularities, via its flow lines. The closed non-singular leaves of this foliation are embedded circles in $\R^3$ and can be studied from a knot theoretical viewpoint. However, not all vector fields have periodic orbits, not even on compact domains. In this case we can still try to study the asymptotical knotting of flow lines. For example, it makes sense to speak of the asymptotic linking number of pairs of orbits of a smooth volume-preserving vector field on a homology sphere~\cite{AK}. Another knot invariant, the signature, was studied by Gambaudo and Ghys in~\cite{GG1}. They proved the existence of an asymptotic signature invariant for orbits of a 
smooth volume-preserving vector field on a compact domain of $\R^3$ and related it to the asymptotic linking number. Recently we could prove the existence of an asymptotic Rasmussen invariant~\cite{B}.
Both the signature and Rasmussen's invariant are so-called concordance invariants. In this note we will study the asymptotics for a family of invariants including many concordance invariants. 

\begin{definition} A link invariant $\tau$ with values in $\R$ is a \emph{linear saddle invariant}, if it satisfies the following two conditions:
\begin{enumerate} 
\item [i)] additivity under disjoint union of links: $\tau(L_1 \sqcup L_2)=\tau(L_1)+\tau(L_2)$, \\
\item [ii)] if two oriented links $L_1$, $L_2$ are related by a saddle point move then: $|\tau(L_1)-\tau(L_2)| \leq C_{\tau}$, where $C_{\tau}>0$ is a constant not depending on $L_1$, $L_2$. \\
\end{enumerate}
\end{definition}

Here a \emph{saddle point move} is a local move that acts on link diagrams as shown in Figure~1.

\begin{figure}[ht]
\scalebox{1}{\raisebox{-0pt}{$\vcenter{\hbox{\epsffile{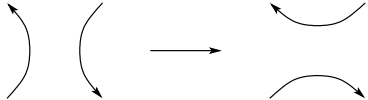}}}$}} 
\caption{saddle point move}
\end{figure} 

Before stating our main result, we have to explain in which sense we want to study asymptotic flow lines: let $X$ be a smooth volume-preserving ergodic vector field on a compact domain $G \subset \R^3$ with smooth boundary. We assume that $X$ is tangent to the boundary $\partial G$ and has only isolated singularities of Morse type, i.e. singularities corresponding to critical points of a Morse function on $\R^3$. The ergodicity of $X$ means that every measurable function on $G$ which is constant on flow lines of $X$ is constant almost everywhere. Let $x \in G$ be a non-periodic point of the flow of $X$, in particular $X(x) \neq 0$. For every $T>0$ there is a piece of flow line starting at $x$ and ending at $\phi_X(T,x)$, where $\phi_X: \R \times G \rightarrow G$ denotes the flow of $X$. We define a subset $K(T,x) \subset \R^3$ as the union of this piece of flow line and the geodesic segment joining $x$ and $\phi_X(T,x)$ in $\R^3$. For almost all $x \in G$, $T>0$, the subset $K(T,x)$ is actually an embedded circle, i.e. a knot in $\R^3$ (see~\cite{GG1}, \cite{V}).

\begin{theorem} Let $\tau$ be a linear saddle invariant with values in $\R$. For almost all $x \in G$ the limit
$$\tau(X,x):=\lim_{T \rightarrow \infty} \frac{1}{T^2} \tau(K(T,x)) \in \R $$
exists.
\end{theorem}

A smooth oriented \emph{cobordism} between two oriented links $L_1$, $L_2 \subset \R^3$ is a smooth oriented surface $S$ relatively embedded in $\R^3 \times [0,1]$ with boundary components $\partial S \cap (\R^3 \times \{0\})=L_1$ and $\partial S \cap (\R^3 \times \{1\})=L_2$. Two oriented links $L_1$, $L_2$ that are locally related by a saddle point move are also related by a smooth oriented cobordism of Euler characteristic $-1$. In fact, this cobordism can even be embedded in $\R^3$ as a saddle surface. Two oriented links are called concordant, if they are related by a cobordism which is a disjoint union of cylinders. A concordance invariant is a link invariant which is constant on equivalence classes of concordant links. Most of the known concordance invariants are linear saddle invariants, for example the signature invariant and Rasmussen's invariant. The latter was first introduced for knots~\cite{Ra} and then for links~\cite{BW}. The classical signature belongs to a family of link invariants called $\omega$-signatures, parametrized by the unit circle in $\C$. The $\omega$-signatures generate an infinite-dimensional vector space of linear saddle invariants (\cite{GG2}, see also~\cite{D} for a calculation of $\omega$-signatures for periodic orbits of the Lorenz flow). Nevertheless, the vector space of asymptotic linear saddle invariants is $1$-dimensional.

\begin{theorem} Let $\tau$ be a linear saddle invariant with values in $\R$. There exists a constant $\alpha \in \R$ such that for almost all $x \in G$:
$$\tau(X,x)=\alpha \sigma(X,x),$$
where $\sigma(X,x)$ is the asymptotic signature invariant.
\end{theorem}

A very special case of Theorem~2 was proved in~\cite{B}: the asymptotic Rasmussen invariant equals twice the asymptotic signature invariant.

The slice genus $g_*(L)$ of an oriented link $L \subset S^3=\partial D^4$ is the minimal genus among all smooth oriented connected surfaces embedded in the $4$-ball $D^4$ with boundary $L$. Unfortunately the slice genus is not a linear saddle invariant, since it is not additive. Nevertheless, there exists an asymptotic slice genus invariant.

\begin{theorem} For almost all $x \in G$ the limit
$$g_*(X,x):=\lim_{T \rightarrow \infty} \frac{1}{T^2} g_*(K(T,x)) \in \R $$
exists and coincides with $|\sigma(X,x)|$.
\end{theorem}

The proofs of Theorems~1 and~2 heavily rely on (and, at the same time, simplify parts of) Gambaudo and Ghys' work, which we will summarize in Section~3. Section~2 contains a lemma on linear saddle invariants that is needed in the proofs of Theorems~1, 2 and~3. In Section~4 we compute the constant $\alpha$ of Theorem~2 for the $\omega$-signatures, where $\omega \in \C$ is a root of unity. Section~5 contains the proof of Theorem~3.

\section{Linear saddle invariants of torus type links}

The signature invariant for links has a good asymptotic behaviour on torus links, in the following sense~\cite{GLM}:

\begin{equation}
\lim_{n \rightarrow \infty} \frac{1}{n^2} \sigma(T(n,n))=\frac{1}{2}.
\label{sigma}
\end{equation}

Here $T(n,n)$ denotes the $n$-component torus link of type $(n,n)$. It is essentially this feature that is responsible for the existence of the asymptotic signature invariant for ergodic vector fields~\cite{GG1}. We will study the behaviour of any linear saddle invariant on a family of links $\{K(n,m)\}$ parametrized by pairs of natural numbers $(n,m)$. For a given pair $(n,m) \in \N \times \N$, we introduce a link $K(n,m)$ as the closure of the following positive braid (see Figure~2, for $n=3$, $m=4$):
$$(\sigma_n \sigma_{n-1} \cdots \sigma_1)(\sigma_{n-1} \sigma_{n-2} \cdots \sigma_2) \cdots (\sigma_{n+m-1} \sigma_{n+m-2} \cdots \sigma_m).$$
Here $\sigma_k$ denotes the $k$-th positive standard generator of the braid group $B_{n+m}$. The link $K(n,m)$ actually coincides with the torus link $T(n,m)$ of type $(n,m)$, but we will not use this fact here.
 
\begin{figure}[ht]
\scalebox{1}{\raisebox{-0pt}{$\vcenter{\hbox{\epsffile{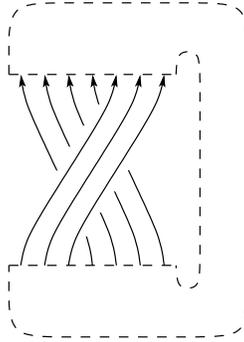}}}$}} 
\caption{the link $K(3,4)$}
\end{figure}

Given any link invariant $\tau$ with values in $\R$, we can define a function $F_{\tau}: \N \times \N \rightarrow \R$ by evaluation on the links $K(n,m)$:
$$F_{\tau}(n,m):=\tau(K(n,m)).$$

\begin{lemma} For every linear saddle invariant $\tau$ with values in $\R$, the limit
$$\lim_{n,m \rightarrow \infty} \frac{1}{nm} F_{\tau}(n,m)=\bar{\tau} \in \R$$
exists, i.e. for all $\epsilon>0$ there exists a natural number $N$, such that $|\frac{1}{nm} F_{\tau}(n,m)-\bar{\tau}| \leq \epsilon$, as soon as $n,m \geq N$.
\end{lemma}

\begin{proof} We may assume that the constant $C_{\tau}$ appearing in condition (ii) is one, since normalization of $\tau$ does not affect the existence of the limit in question. Further, by replacing $\tau(L)$ by $\tau(L)-\# L \, \tau(O)$, we may assume that $\tau(O)=0$. Here $\#L$ denotes the number of components of the link $L$ and $O$ denotes the trivial knot. We then observe that $|F_{\tau}(n,m)| \leq nm$. Indeed, the link $K(n,m)$ has a diagram with $nm$ crossings, and any connected link diagram can be transformed into a trivial knot by applying a saddle point move at some of its crossings, as shown in Figure~3.

\begin{figure}[ht]
\scalebox{1}{\raisebox{-0pt}{$\vcenter{\hbox{\epsffile{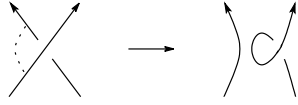}}}$}} 
\caption{deleting a crossing by a saddle point move}
\end{figure}

The function $f_{\tau}: \N \times \N \rightarrow \R$, defined by
$$f_{\tau}(n,m):=\frac{1}{nm} F_{\tau}(n,m),$$
is therefore absolutely bounded: $|f| \leq 1$. We have to show that $f$ has a limit.

Let $n_1,n_2,m \in \N$. There is a sequence of $m$ saddle point moves between the link $K(n_1+n_2,m)$ and the disjoint union of the two links $K(n_1,m)$, $K(n_2,m)$, see Figure~4, for $n_1=3$, $n_2=2$, $m=4$ (the dashed line indicates the area where $4$ successive saddle point moves have to be performed). Therefore:
$$|F_{\tau}(n_1+n_2,m)-F_{\tau}(n_1,m)-F_{\tau}(n_2,m)| \leq m,$$
i.e. the restriction of $F$ to one variable is a quasi-morphism on $\N$ (see \cite{AC} for a detailed account on quasi-morphisms on $\Z$). Applying the above estimate $p-1$ times, we obtain:
$$|F_{\tau}(pn,m)-pF_{\tau}(n,m)| \leq (p-1)m \leq pm,$$
$$|f_{\tau}(pn,m)-f_{\tau}(n,m)| \leq \frac{1}{n}.$$
An analogous estimate holds for the second variable:
$$|f_{\tau}(n,qm)-f_{\tau}(n,m)| \leq \frac{1}{m}.$$
The last two inequalities imply the existence of a limit for $f$: let $\epsilon>0$, $N \in \N$, $\frac{4}{N} \leq \epsilon$. Then, for all $p,q,n,m \in \N$, $p,q,n,m \geq N$: 
\begin{align*}
|f_{\tau}(p,q)-f_{\tau}(n,m)| 
& \leq |f_{\tau}(n,m)-f_{\tau}(pn,qm)|+|f_{\tau}(pn,qm)-f_{\tau}(p,q)| \\
& \leq \frac{1}{n}+\frac{1}{m}+\frac{1}{p}+\frac{1}{q} \leq \frac{4}{N} \leq \epsilon
\end{align*}

\begin{figure}[ht]
\scalebox{1}{\raisebox{-0pt}{$\vcenter{\hbox{\epsffile{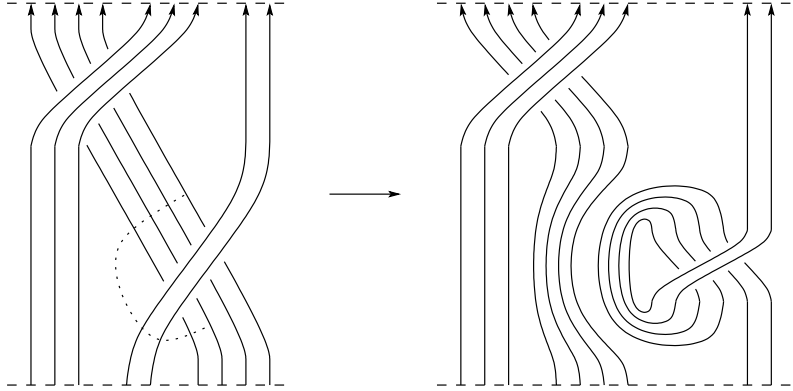}}}$}} 
\caption{the links $K(5,4)$ and $K(3,4) \sqcup K(2,4)$}
\end{figure}

\end{proof}

\begin{remark} 
We could define the family of links $\{K(n,m)\}$ as closures of negative braids, as well. The limit of Lemma~1 would thereby change its sign. Indeed, the disjoint union of the link $K(n,m)$ and its mirror image can be transformed into a trivial link with $n+m$ components by a sequence of $n+m$ saddle point moves, hence into a trivial knot by $2n+2m-1$ saddle point moves (this is an easy exercise). Thus the values of any linear saddle invariant on $K(n,m)$ and its mirror image differ by their signs, up to an affine error in $n$ and $m$. 
\end{remark}

\section{Proofs of Theorems~1 and~2}

Gambaudo and Ghys' proof of the existence of an asymptotic signature invariant is actually well-adapted for any linear saddle invariant. Under the assumptions on the vector field $X$, they cover the complement of the singularities of $X$ (finitely many in number) by a countable collection of flow boxes $\{\mathcal{F}_i\}_{i \in \N}$. An important feature of this collection is that the flow time of each box, i.e. the time it takes to pass through the box, is bounded from below by a global constant $\lambda>0$. For all $x \in G$, $T>0$, they define a number $n_i(T,x)$ which measures how many times the flow line starting at $x$ and ending at $\phi_X(T,x)$ enters the flow box $\mathcal{F}_i$. The following estimates are obvious:
$$0 \leq n_i(T,x) \leq \frac{T}{\lambda}.$$
Using Birkhoff's ergodic theorem, Gambaudo and Ghys argue that for almost all $x \in G$ the limit
\begin{equation}
\lim_{T \rightarrow \infty} \frac{1}{T} n_i(T,x)=n_i(x)>0
\label{averagenumber}
\end{equation}
exists (and is proportional to the volume of the flow box $\mathcal{F}_i$). 

The family of flow boxes $\{\mathcal{F}_i\}$ comes together with a good projection $\pi_0: \R^3 \rightarrow \R^2$ onto a plane. For every $\epsilon>0$, there exists a finite subset $\mathcal{C} \subset \N$, such that for almost all $x \in G$, $T>0$ large enough, the diagram $\pi_0(K(T,x))$ is regular and has a `large' and a `small' part: up to an error of $\epsilon T^2$, the crossings of $\pi_0(K(T,x))$ arise from pairs of overcrossing flow boxes $\mathcal{F}_i$, $\mathcal{F}_j$, for $i,j \in \mathcal{C}$. At the spots of overcrossing flow boxes (again finitely many in number, say $c_1, \ldots, c_N$), the diagram $\pi_0(K(T,x))$ looks locally like a rectangular grid with $n_i(T,x) n_j(T,x)$ crossings, see the first diagram of Figure~5. 

The crucial step that allows us to prove Theorem~1 is to split off links of type $\widetilde{K}(n_i(T,x),n_j(T,x))$ at every crossing $c_k$ ($1 \leq k \leq N$), by a controlled number of saddle point moves. Here $i$, $j \in \mathcal{C}$ denote the indices of the flow boxes crossing at $c_k$, and $\widetilde{K}(n_i(T,x),n_j(T,x))$ denotes either a link of type $K(n_i(T,x),n_j(T,x))$ or a mirror image thereof. This can be done by applying precisely $n_i(T,x))+n_j(T,x)$ saddle point moves at $c_k$, as illustrated in Figure~5. Here again, the dashed line on the left (resp. on the right) indicates the area where $n_i$ (resp. $n_j$) successive saddle point moves have to be performed.

\begin{figure}[ht]
\scalebox{0.8}{\raisebox{-0pt}{$\vcenter{\hbox{\epsffile{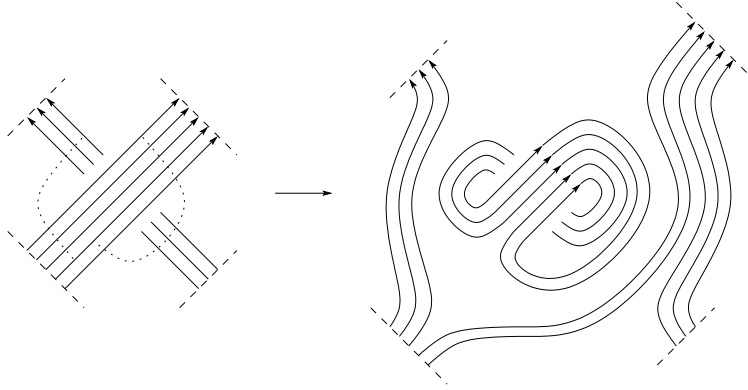}}}$}} 
\caption{splitting off $K(n_i,n_j)$}
\end{figure}

Altogether we need at most $N(\frac{T}{\lambda}+\frac{T}{\lambda})=\frac{2N}{\lambda}T$ saddle point moves to transform the knot $K(T,x)$ into a disjoint union of two links $L_1(T,x) \sqcup L_2(T,x)$, where $L_1(T,x)$ is a disjoint union of $N$ links of type $\widetilde{K}(n_i,n_j)$ and $L_2(T,x)$ is a link whose crossing number is at most $\epsilon T^2$. Since every saddle point move produces at most one new link component, the link $L_2(T,x)$ has at most $\frac{2N}{\lambda}T$ components. Therefore, it can be transformed into a trivial knot by $\frac{2N}{\lambda}T+\epsilon T^2$ saddle point moves, at most. The following two estimates hold for any linear saddle invariant $\tau$ with normalization $C_{\tau}=1$ and $\tau(O)=0$:
$$|\tau(K(T,x))-\tau(L_1(T,x) \sqcup L_2(T,x))| \leq \frac{2N}{\lambda}T,$$
$$|\tau(L_1(T,x) \sqcup L_2(T,x))-\tau(L_1(T,x))|=|\tau(L_2(T,x)| \leq \frac{2N}{\lambda}T+\epsilon T^2,$$
These two inequalities together with the equality
$$\tau(L_1(T,x))=\sum_{k=1}^{N} \tau(\widetilde{K}(n_i(T,x),n_j(T,x)))$$ imply:
$$|\frac{1}{T^2}\tau(K(T,x))-\frac{1}{T^2}\sum_{k=1}^{N} \tau(\widetilde{K}(n_i(T,x),n_j(T,x)))| \leq \frac{4N}{\lambda T}+\epsilon.$$
Therefore it remains to show that the limit 
$$\lim_{T \rightarrow \infty} \frac{1}{T^2} \tau(\widetilde{K}(n_i(T,x),n_j(T,x)))$$
exists, for all $i$, $j \in \mathcal{C}$, for almost all $x \in G$. This is an easy consequence of Lemma~1, the remark of the same section, and (\ref{averagenumber}). The limit turns out to be $\pm n_i(x) n_j(x) \bar{\tau}$.

As the proof of Theorem~1 shows, two linear saddle invariants $\tau_1$, $\tau_2$ with
$$\lim_{n,m \rightarrow \infty} \frac{1}{nm} F_{\tau_1}(n,m)=\alpha \lim_{n,m \rightarrow \infty} \frac{1}{nm} F_{\tau_2}(n,m),$$
for some $\alpha \in \R$, lead to proportional asymptotic invariants: $\tau_1(X,x)=\alpha \tau_2(X,x)$, for almost all $x \in G$. Therefore the vector space of asymptotic linear saddle invariants is $1$-dimensional, generated by the asymptotic signature, as stated in Theorem~2.

\section{Asymptotic $\omega$-signatures}

The classical signature invariant $\sigma$ of a link $L$ is defined as the signature of any symmetrized Seifert matrix $V$ for $L$:
$$\sigma(L)=\text{sign} (V+V^T).$$
More generally, there exists a link invariant $\sigma_{\omega}$ for every $\omega \in \C$, $|\omega|=1$. It is defined as the number of positive eigenvalues minus the number of negative eigenvalues of the hermitian matrix
$$(1-\omega)V+(1-\bar{\omega})V^T.$$
All these invariants are easily seen to be linear saddle invariants. By Theorem~2, the corresponding asymptotic invariants are multiples of the asymptotic signature invariant. We will determine the explicit ratio $\alpha \in \R$ in case $\omega$ is a root of unity.

\begin{proposition} Let $X$ be a smooth vector field on a compact domain $G \subset \R^3$, as in the introduction, and let $\omega=e^{2 \pi i \theta} \in \C$ be a root of unity. For almost all $x \in G$
$$\sigma_{\omega}(X,x)=(4 \theta (1-\theta)) \sigma(X,x).$$
\end{proposition}

\begin{proof}
By the proof of Theorem~1 we have to show
$$\lim_{n \rightarrow \infty} \frac{1}{n^2} \sigma_{\omega}(K(n,n))=(4 \theta (1-\theta)) \lim_{n \rightarrow \infty} \frac{1}{n^2} \sigma(K(n,n)).$$
This follows from a calculation of $\omega$-signatures for torus links which was carried out by Gambaudo and Ghys \cite{GG2}. The torus link of type $(n,p)$, $n \geq 2$, $p \geq 1$, is defined as the closure of the positive braid $(\sigma_1 \sigma_2 \cdots \sigma_{n-1})^p \in B_n$. Let $\omega=e^{2 \pi i \theta}$, where $\theta \in \Q$, $0 \leq \theta<1$, and let $l_{\theta} \in \{1,2,\ldots,n\}$ be the unique natural number with $\frac{l_{\theta}-1}{n} \leq \theta < \frac{l_{\theta}}{n}$. Proposition 5.2 of \cite{GG2} immediately implies
$$|\sigma_{\omega}(T(n,p))-2p\theta (n+1-2l_{\theta})-\frac{2p}{n} l_{\theta}(l_{\theta}-1)| \leq 2n,$$
\begin{equation}
\lim_{n \rightarrow \infty} \frac{1}{n^2} \sigma_{\omega}(T(n,n))=2\theta-4\theta^2+2\theta^2=2\theta(1-\theta).
\label{toruslimit}
\end{equation}
As we mentioned in Section~2, the two families of links $\{K(n,m)\}$ and $\{T(n,m)\}$ coincide. In particular, the links $K(n,n)$ and $T(n,n)$ are isotopic, for all $n \geq 1$. Indeed, both links are isotopic to the closure of $n$ parallel strands with a full twist, see Figure~6, for $n=4$. Thus the equation (\ref{toruslimit}) holds for $K(n,n)$, as well:
\begin{equation}
\lim_{n \rightarrow \infty} \frac{1}{n^2} \sigma_{\omega}(K(n,n))=2\theta(1-\theta).
\label{Klimit}
\end{equation}
For $\theta=\frac{1}{2}$, i.e. for the classical signature $\sigma=\sigma_{-1}$, this limit equals $\frac{1}{2}$, in accordance with (\ref{sigma}). We conclude
$$\lim_{n \rightarrow \infty} \frac{1}{n^2} \sigma_{\omega}(K(n,n))=(4 \theta (1-\theta)) \lim_{n \rightarrow \infty} \frac{1}{n^2} \sigma(K(n,n)).$$
\end{proof}

\begin{figure}[ht]
\scalebox{1}{\raisebox{-0pt}{$\vcenter{\hbox{\epsffile{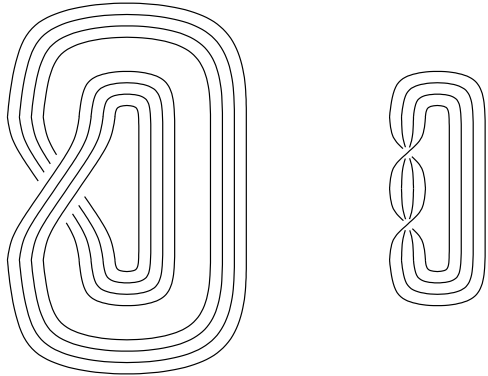}}}$}} 
\caption{the isotopic links $K(4,4)$ and $T(4,4)$}
\end{figure}

\section{Asymptotic slice genus}

The slice genus of links is difficult to determine, in general. In the case of positive braid links, the slice genus was first determined by Kronheimer and Mrowka~\cite{KM} (see also~\cite{Ru}). An alternative method was recently given by Rasmussen~\cite{Ra}. Their results imply the existence of the following limit:
\begin{equation}
\lim_{n,m \rightarrow \infty} \frac{1}{nm} g_*(K(n,m))=\frac{1}{2}.
\label{genuslimit}
\end{equation}
Thus Lemma~1 is true for the slice genus. However, in contrast to linear saddle invariants, the slice genus is invariant under mirror image. Therefore the limit (\ref{genuslimit}) stays the same if we replace the links $K(n,m)$ by their mirror images. The fact that the slice genus is not additive under disjoint union of links forces us to adapt the proof of Theorem~1. Thereby we are still allowed to use the second condition of linear saddle invariants: if two oriented links $L_1$, $L_2$ are related by a saddle point move then
\begin{equation}
|g_*(L_1)-g_*(L_2)| \leq 1.
\label{genusestimate}
\end{equation}
In the following we keep the notation of Section~3.

First we transform the knot $K(T,x)$ into a disjoint union of two links $L_1(T,x) \sqcup L_2(T,x)$, by a sequence of $z \leq \frac{2N}{\lambda}T$ saddle point moves. The link $L_1(T,x) \sqcup L_2(T,x)$ in turn can be transformed into the link $L_1(T,x)$ by a sequence of $z \leq \frac{2N}{\lambda}T+\epsilon T^2$ saddle point moves. Thus
$$|g_*(K(T,x))-g_*(L_1(T,x))| \leq \frac{4N}{\lambda}T+\epsilon T^2.$$
Now comes the place where we run into trouble: we cannot apply the additivity of $g_*$ to the link $L_1(T,x)$, which is a disjoint union of $N$ links of type $\widetilde{K}(n_i(T,x),n_j(T,x))$. Here we use the convention that $\widetilde{K}$ denotes either the link $K$ or its mirror image. We need the following lemma to go around the additivity.

\begin{lemma} Let $a$, $n$, $m$ be natural numbers, $1 \leq a \leq n,m$. There exists a natural number $b$, such that the link $\widetilde{K}(n,m)$ can be transformed into the link $\widetilde{K}(a,b)$ by a sequence of $z \leq m+n+\frac{mn}{a}+am$ saddle point moves.
\end{lemma}

\begin{proof} There exist (unique) natural numbers $k$, $r$ with $r<a$ and $n=ak+r$. In the following, the symbol
$$K \stackrel{x}{\longrightarrow} L$$
means that the link $K$ can be transformed into the link $L$ by a sequence of $x$ saddle point moves, at most. The proof of Lemma~1 implies
$$\widetilde{K}(n,m) \stackrel{m}{\longrightarrow} \widetilde{K}(ak,m) \sqcup \widetilde{K}(r,m) \stackrel{y}{\longrightarrow} \widetilde{K}(a,km) \sqcup \widetilde{K}(r,m),$$
where $y=(k-1)m+(k-1)a$. Further we evidently have
$$\widetilde{K}(a,km) \sqcup \widetilde{K}(r,m) \stackrel{rm}{\longrightarrow} \widetilde{K}(a,km).$$
Altogether these arrows imply
$$\widetilde{K}(n,m) \stackrel{z}{\longrightarrow} \widetilde{K}(a,b),$$
where $b=km$ and
\begin{align*}
z&=m+(k-1)(m+a)+rm\\
& \leq m+\frac{n}{a}(m+a)+am\\
& \leq m+n+\frac{mn}{a}+am\\
\end{align*}
\end{proof}

We will apply Lemma~2 to the links $\widetilde{K}(n_i(T,x),n_j(T,x))$ and suitable natural numbers $a(T,x)$. Let $l \in \mathcal{C}$ be the index for which $n_l(x)$ is minimal (here again $\mathcal{C} \subset \N$ and $n_i(x) \in \R$ are defined as in Section~3).
Let $a(T,x)$ be a natural number `close' to $\sqrt{n_l(T,x)}$, for example the integral part of $\sqrt{n_l(T,x)}$. Further let $\widetilde{K}(n_i(T,x),n_j(T,x))$ be any link component of $L_1(T,x)$. According to Lemma~2, there exists a natural number $b_{ij}(T,x)$ with
$$\widetilde{K}(n_i(T,x),n_j(T,x)) \stackrel{z}{\longrightarrow} \widetilde{K}(a(T,x),b_{ij}(T,x)),$$
where $z \leq n_i+n_j+\frac{n_i n_j}{a}+an_j$ (here we suppress the parameters $T$ and $x$). By our construction, the numbers $a(T,x)$ grow like $\sqrt{T}$. Combining this with (\ref{genusestimate}), we obtain
$$\lim_{T \rightarrow \infty} \frac{1}{T^2} |g_*(\widetilde{K}(n_i(T,x),n_j(T,x)))-g_*(\widetilde{K}(a(T,x),b_{ij}(T,x)))|=0.$$
Thus we can replace the link $L_1(T,x)$ by a disjoint union of $N$ links of type $\widetilde{K}(a(T,x),b_{ij}(T,x))$ with the same $a(T,x)$ for all $N$ components. This union can further be transformed into one single link of type 
$\widetilde{K}(a(T,x),b(T,x))$ by $z \leq N a(T,x)$ saddle point moves. Here $b(T,x)$ is a sum with signs of all $b_{ij}(T,x))$. The numbers $b_{ij}(T,x)$, $b(T,x)$ grow like $T \sqrt{T}$. Therefore the limit
$$\lim_{T \rightarrow \infty} \frac{1}{T^2} a(T,x)b(T,x)$$
exists. From this and (\ref{genuslimit}) we deduce the existence of the limit
$$\lim_{T \rightarrow \infty} \frac{1}{T^2} g_*(K(T,x))=\lim_{T \rightarrow \infty} \frac{1}{T^2} g_*(\widetilde{K}(a(T,x),b(T,x))).$$
The last statement of Theorem~3 is obvious since
$$\lim_{n,m \rightarrow \infty} \frac{1}{nm} g_*(\widetilde{K}(n,m))=\lim_{n,m \rightarrow \infty} \frac{1}{nm} |\sigma(\widetilde{K}(n,m))|=\frac{1}{2}.$$

\textbf{Acknowledgements}. I would like to thank David Cimasoni for teaching me concordance invariants and for finding an adequate definition of linear saddle invariants.

\bigskip
\noindent
Department of Mathematics,
ETH Z\"urich, 
Switzerland

\bigskip
\noindent
\emph{sebastian.baader@math.ethz.ch}

\end{document}